\documentclass[aps,pre,reprint,amsmath,amssymb]{revtex4-2}
\usepackage{txfonts}
\usepackage{algorithmic}
\usepackage{algorithm}
\usepackage{graphicx}
\usepackage{bm}
\usepackage{url}

\newcommand{\argmin}{\mathop{\rm arg~min}\limits}

\begin{document}

\title{Improvement of system identification of stochastic systems via Koopman generator and locally weighted expectation}

\author{Yuki Tahara, Kakutaro Fukushi, Shunta Takahashi, Kayo Kinjo, and Jun Ohkubo}

\affiliation{Graduate School of Science and Engineering, Saitama University, Sakura, Saitama, 338--8570 Japan}

\begin{abstract}
The estimation of equations from data is of interest in physics. One of the famous methods is the sparse identification of nonlinear dynamics (SINDy), which utilizes sparse estimation techniques to estimate equations from data. Recently, a method based on the Koopman operator has been developed; the generator extended dynamic mode decomposition (gEDMD) estimates a time evolution generator of dynamical and stochastic systems. However, a naive application of the gEDMD algorithm cannot work well for stochastic differential equations because of the noise effects in the data. Hence, the estimation based on conditional expectation values, in which we approximate the first and second derivatives on each coordinate, is practical. A naive approach is the usage of locally weighted expectations. We show that the naive locally weighted expectation is insufficient because of the nonlinear behavior of the underlying system. For improvement, we apply the clustering method in two ways; one is to reduce the effective number of data, and the other is to capture local information more accurately. We demonstrate the improvement of the proposed method for the double-well potential system with state-dependent noise.
\end{abstract}

\maketitle

\section{Introduction}
\label{introduction}

Recent developments in machine learning enable us to estimate system equations from data. The estimation is easy for linear systems but generally difficult for nonlinear systems. One notable method is the sparse identification of nonlinear dynamics (SINDy) \cite{Brunton2016}, which utilizes sparse estimation techniques to infer equations from data for deterministic dynamical systems. Some extensions of the SINDy to stochastic systems have already been carried out \cite{Boninsegna2018,Huang2022,Wanner2024}.

Here, we focus on the method based on the Koopman operator method. The Koopman operator \cite{Koopman1931} is applicable in various research fields, such as molecular dynamics, meteorology, and economics. With the aid of the Koopman operator, it becomes possible to analyze nonlinear systems using linear techniques. However, we cannot use it on computers since the Koopman operator resides in an infinite-dimensional space. Therefore, it is necessary to approximate the Koopman operator in finite dimensions. For this purpose, some algorithms have been proposed; the dynamic mode decomposition (DMD) \cite{Schmid2010} extracts information about the Koopman operator from the time-series data of the system. The mode decomposition is also available for a partial differential equation, such as the nonlinear Burgers equation \cite{Nakao2020}. The extended dynamic mode decomposition (EDMD) \cite{Williams2015} is an extension of DMD, which utilizes a dictionary of basis functions to approximate the Koopman operator in finite dimensions. Applications of the EDMD include the prediction of the system's future states, system identification \cite{Mauroy2020}, and control \cite{Kaiser2017}. There are many research papers on these topics; see Ref.~\cite{Brunton2022} for more details.

Moreover, an extension of EDMD called generator EDMD (gEDMD) \cite{Klus2020} is applicable to estimate stochastic differential equations from data. As described in Ref.~\cite{Klus2020}, numerical experiments were performed for a double-well potential system, which recovered the correct underlying stochastic differential equation. However, as we will explain later, the dataset used in Ref.~\cite{Klus2020} was not simple sample trajectories. When we apply the gEDMD algorithm under the assumption of trajectory data only, the noise in the data considerably affects the estimation; we will discuss this in a numerical demonstration. The Koopman operator and the related practical algorithms have attracted attention as a new analytical method for nonlinear systems, and improving the gEDMD method is highly significant for other future research.

In this study, we aim to improve the accuracy of the gEDMD estimation. Hence, we propose some ideas of preprocessing for noisy time-series data to mitigate the effect of noise. One of them is the practical implementation of the calculation of conditional expectations as weighted expectations with kernel functions. We also apply clustering procedures of data in two ways. One is to reduce the computational cost; the other is to capture the local features of data adequately. We demonstrate the proposed methods with the double-well potential system with the state-dependent noise used in the previous study \cite{Klus2020}.

The structure of this paper is as follows. In Sect.~2, we describe the problem settings and a concrete model; the concrete example will help readers understand the following discussions. In Sect.~3, we briefly review the gEDMD algorithm. In Sect.~4, we present our proposals, i.e., weighted expectation using kernel functions and clustering-based data preprocessing. In Sect.~5, we numerically demonstrate the proposed methods. Finally, in Sect.~6, we provide the concluding remarks and discussions for future work.

\section{Problem Settings and Model}

Here, we briefly explain the problem settings. In this study, we employed almost the same settings and the model as in Ref.~\cite{Klus2020}. In the following sections, we sometimes used the double-well potential example, which helps readers understand the explanations and discussions.

\subsection{Stochastic differential equation}

We here focus on the following type of stochastic differential equation:
\begin{align}
d\bm{X}(t) = \bm{b}(\bm{X}(t))dt + \Sigma(\bm{X}(t)) d\bm{W}(t),
\label{eq_sde}
\end{align}
where $\bm{b}(\bm{x}): \mathbb{R}^{D}\rightarrow\mathbb{R}^{D}$ is a $D$-dimensional vector function, which yields the drift coefficients, and $\Sigma(\bm{x}): \mathbb{R}^{D}\rightarrow\mathbb{R}^{D\times s}$ is a matrix-valued function for diffusion coefficients. $\bm{W}(t)$ is an $s$-dimensional Wiener process. Note that the drift and diffusion coefficient functions are nonlinear in general. For more details of stochastic differential equations, see, for example, Ref.~\cite{Gardiner_book}. The aim of this work is to estimate the drift and diffusion coefficient functions only from trajectory data of $\bm{X}(t)$.

\subsection{Double-well potential}

\begin{figure}
    \begin{center}
        \includegraphics[width=80mm]{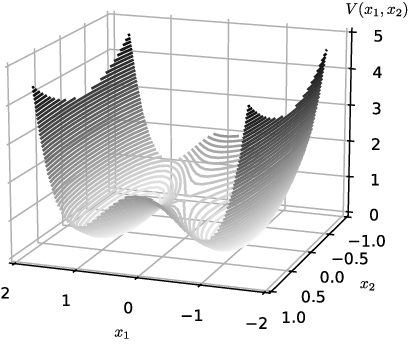}
        \caption{Potential of the problem. There is a double-well structure. As we see later, the actual distribution deviates slightly from the potential contours owing to state-dependent noise.}
        \label{fig_potential}
    \end{center}
\end{figure}

In Ref.~\cite{Klus2020}, the following two-variable system was used to demonstrate the gEDMD algorithm:
\begin{align}
V(\bm{x})=(x_1^2-1)^2+x_2^2,
\label{eq_potential}
\end{align}
which leads to the drift coefficient $\bm{b}(\bm{x})$ in Eq.~\eqref{eq_sde} as follows:
\begin{align}
\bm{b}(x)=-\nabla{V}(x)=\begin{bmatrix}
4x_1-4{x_1}^3 \\
-2x_2 \\
\end{bmatrix}.
\label{eq_drift}
\end{align}
Figure~\ref{fig_potential} illustrates the double-well potential. There are two wells, and $x_1$ and $x_2$ will move back and forth between the two wells over time. However, there is a state-dependent noise as follows:
\begin{align}
\Sigma(\bm{x})=\begin{bmatrix}
0.7 & x_1 \\
0 & 0.5 \\
\end{bmatrix}.
\end{align}
Thus, we expect that the actual distribution deviates slightly from the potential contours.

\begin{figure}
    \begin{center}
        \includegraphics[width=75mm]{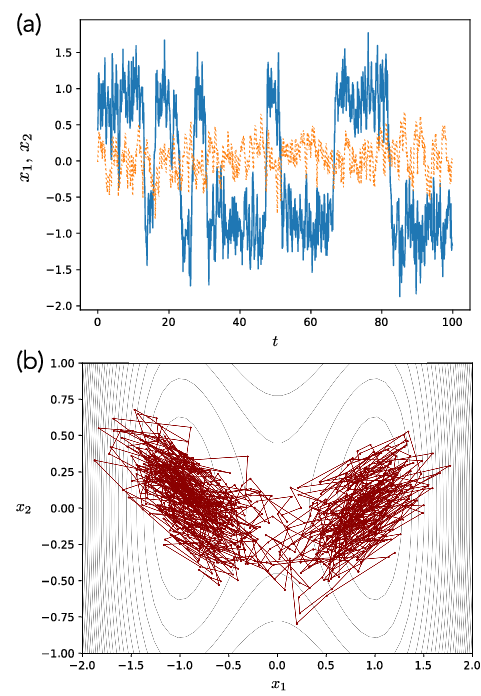}
        \caption{(Color online) A sample trajectory from the double-well potential model. (a) Time-domain picture. The solid line corresponds to $x_1$, and the dashed line corresponds to $x_2$. (b) A sample trajectory depicted on the state space. The contour lines represent the double-well potential.}
        \label{fig_sample}
    \end{center}
\end{figure}

Figure~\ref{fig_sample}(a) shows the time evolution of $x_1$ and $x_2$ for the double-well potential with the state-dependent diffusion term. We clearly see that the time evolution of $x_1$ shows the move between the two wells around $x_1=1.0$ and $x_1=-1.0$. Figure~\ref{fig_sample}(b) illustrates the sample trajectory on the state space. The contour lines represent the magnitude of the potential. Note that the diffusion term is state-dependent, and therefore the trajectory does not follow the contour lines but is oriented diagonally.

\section{Equation Estimation via the Koopman Generator}

As for the details of the Koopman operator, please see, for example, Ref.~\cite{Brunton2022}. In this section, we will briefly look at only the minimum required for equation estimation.

\subsection{Koopman operator and generator}

The Koopman operator acts on a function, not a state vector, and it yields a function that gives an observed quantity after time evolution. Let $f(\bm{x}) \in \mathcal{F}$ be an observable function, where $\mathcal{F}$ is a function space for observables. Then, the Koopman operator $\mathcal{K}^{\Delta t}$ acts on the function $f(\bm{x})$ as 
\begin{align}
\left(\mathcal{K}^{\Delta t} f \right) (\bm{x})
= \mathbb{E} \left[ \left. f(\bm{X}(t+\Delta t)) \, \right| \, \bm{X}(t)=\bm{x}\right],
\end{align}
where $\mathbb{E}[ \cdot ]$ represents the expectation value. Note that the state before the time evolution is $\bm{X}(t)=\bm{x}$. The new function $\mathcal{K}^{\Delta t} f$ takes the state before time evolution as its argument, and it yields the expectation value of function $f$ after the time evolution with $\Delta t$. Hence, the Koopman operator gives the time evolution in the function space instead of the state space.

In the gEDMD algorithm \cite{Klus2020}, we seek the operator that gives the time derivative of $\mathcal{K}^{\Delta t} f$. Then, we introduce the Koopman generator $\mathcal{L}: \mathcal{F}\rightarrow\mathcal{F}$, which acts on function $f$ to produce another function, i.e., the time derivative. The definition of the Koopman generator $\mathcal{L}$ associated with the Koopman operator $\mathcal{K}^{\Delta t}$ is given as 
\begin{align}
(\mathcal{L}f)(\bm{x})&=\lim_{\Delta t\to0}\frac{(\mathcal{K}^{\Delta t}f)(\bm{x})-f(\bm{x})}{\Delta t}.
\end{align}

When we consider Eq.~\eqref{eq_sde} and a second-order continuously differentiable observable $f$, the Koopman generator $\mathcal{L}$ acts as \cite{Klus2020}:
\begin{align}
(\mathcal{L}f)(\bm{x})
&=\sum_{d=1}^{D} b_d(\bm{x})\frac{\partial}{\partial x_d}f(\bm{x})+\frac{1}{2}\sum_{i=1}^{D}\sum_{j=1}^{D}a_{ij}(\bm{x})\frac{\partial^2}{\partial x_i\partial x_j}f(\bm{x}),
\label{eq_generator}
\end{align}
where $A(\bm{x}) = [a_{ij}(\bm{x})]$ is a $D\times D$ dimensional square matrix defined by $A(\bm{x})=\Sigma(\bm{x})\Sigma(\bm{x})^{\top}$. Note that the adjoint operator of $\mathcal{L}$ is 
\begin{align}
\mathcal{L}^{*}
&=\sum_{d=1}^{D} \frac{\partial}{\partial x_d} b_d(\bm{x}) 
+\frac{1}{2}\sum_{i=1}^{D}\sum_{j=1}^{D}\frac{\partial^2}{\partial x_i\partial x_j}  a_{ij}(\bm{x}).
\label{eq_fokker_planck}
\end{align}
We see that the operator $\mathcal{L}^{*}$ yields the Fokker--Planck equation to describe the time evolution of the probability density function. Actually, the operator in Eq.~\eqref{eq_generator} corresponds to that in the Kolmogorov backward equation \cite{Gardiner_book}.

\subsection{Derivation of the Koopman generator matrix $L$}

Assume that we have $N$ data points $\{\bm{x}_{n}\}_{n=1}^{N}$ in time-series data. Our aim here is to estimate $\bm{b}(\bm{x})$ and $A(\bm{x})$ in Eq.~\eqref{eq_generator}. Then, we can estimate the diffusion coefficient $\Sigma(\bm{x})$ in Eq.~\eqref{eq_sde} via the Cholesky decomposition of $A(\bm{x})$, although the decomposition is not unique.

We employ a set of basis functions to handle the observable function $f(\bm{x})$. The set is called a dictionary, and we denote it as a vector function $\bm{\psi}(\bm{x})$. For simplicity, we here employ the monomial dictionary. The double-well potential example has two state variables and the following dictionary is available:
\begin{align}
&\bm{\psi}(\bm{x}) \nonumber \\
&= [\psi_1(\bm{x}), \psi_2(\bm{x}),\psi_3(\bm{x}), \psi_4(\bm{x}), 
\psi_5(\bm{x}),\psi_6(\bm{x}),\psi_7(\bm{x}),\psi_8(\bm{x}), \dots]^{\top}\nonumber \\
&= [1, \, x_1, \, x_2, \, x_1^2, \, x_1x_2, \, x_2^2, \, x_1^3, \, x_1^2x_2, \, \dots]^{\top}.
\label{eq_dictionary_example}
\end{align}
It is also possible to include other types of functions such as $\sin$ and $\cos$. Of course, we cannot use an infinite number of functions in practical computation, and some approximation is necessary. We use a finite number of functions as the dictionary. We denote the number of dictionary functions as $N_{\mathrm{dic}}$.

Next, we introduce the following notation for each data point:
\begin{align}
d\psi_k(\bm{x}_n)
=&
\sum_{d=1}^{D} b_d(\bm{x}_n)\frac{\partial}{\partial x_d}\psi_k(\bm{x}_n) \nonumber \\
& +
\frac{1}{2}\sum_{i=1}^{D}\sum_{j=1}^{D}a_{ij}(\bm{x}_n)\frac{\partial^2}{\partial x_i\partial x_j}\psi_k(\bm{x}_n).
\label{eq_dpsi}
\end{align}
Note that the concrete data point $\bm{x}_n$ is used, and hence Eq.~\eqref{eq_dpsi} yields a scalar value. Then, we construct the following data vector for data $n$:
\begin{align}
\bm{\psi}(\bm{x}_n) =
\begin{bmatrix}
\psi_1(\bm{x}_n) \\
\vdots \\
\psi_{N_{\mathrm{dic}}}(\bm{x}_n)
\end{bmatrix}, \quad
\bm{d\psi}(\bm{x}_n) =
\begin{bmatrix}
d\psi_1(\bm{x}_n) \\
\vdots \\
d\psi_{N_{\mathrm{dic}}}(\bm{x}_n)
\end{bmatrix}.
\end{align}

The Koopman generator $\mathcal{L}$ is infinite-dimensional since it acts on functions. Hence, it is infeasible to use it on a computer. Therefore, we need an approximation of the Koopman generator $\mathcal{L}$ with a finite-dimensional matrix. Let $L\in\mathbb{R}^{N_{\mathrm{dic}}\times N_{\mathrm{dic}}}$ be the matrix representation of the Koopman generator $\mathcal{L}$ in terms of the specific dictionary $\bm{\psi}$, and consider the following least squares problem:
\begin{align}
L&=\argmin_{\widetilde{L}}\sum_{n=1}^{N}\left\|\bm{d\psi}(\bm{x}_n)-\widetilde{L}\bm{\psi}(\bm{x}_n)\right\|^2.
\label{eq_cost}
\end{align}
The solution to this least squares problem is given as
\begin{align}
L=FG^+,
\label{eq_L}
\end{align}
where
\begin{align}
F=\frac{1}{N}\sum^{N}_{n=1}\big(\bm{d\psi}(\bm{x}_n)\big)\bm{\psi}(\bm{x}_n)^{\top}, \,\,
G=\frac{1}{N}\sum^{N}_{n=1}\bm{\psi}(\bm{x}_n)\bm{\psi}(\bm{x}_n)^{\top},
\end{align}
and $G^+$ denotes the pseudo-inverse of the matrix $G$. We call $L$ the Koopman generator matrix.

There is a comment on the difference in the definition of the Koopman generator from that in Ref.~\cite{Klus2020}. In Ref.~\cite{Klus2020}, the transposed matrix of $L$ was mainly discussed. Although it is natural to employ the above definitions in this  paper, it would be better to consider the transposed matrix $L^{\top}$ for comparison with Ref.~\cite{Klus2020}. Hence, in the later numerical experiments, we mainly discuss $L^{\top}$.

\subsection{Estimation of drift and diffusion coefficient functions}

In the EDMD and gEDMD algorithms, the vector function $\bm{g}(\bm{x}) = \bm{x}$ is called the full-state observable, which directly corresponds to the state variables. For the example of the double-well potential, $\psi_2(\bm{x}) = x_1$ and $\psi_3(\bm{x}) = x_2$ in Eq.~\eqref{eq_dictionary_example} correspond to the full-state observable; i.e., we define
\begin{align}
\bm{g}(\bm{x})=
\begin{bmatrix}
\psi_2(\bm{x})\\
\psi_3(\bm{x}) 
\end{bmatrix}
=
\begin{bmatrix}
x_1 \\
x_2
\end{bmatrix}.
\end{align}

Using the full-state observable, we easily obtain the drift coefficient function $\bm{b}(\bm{x})$ as
\begin{align}
(\mathcal{L}\bm{g})(\bm{x}) = \bm{b}(\bm{x}),
\end{align}
i.e., the action of the Koopman generator $\mathcal{L}$ to the full-state observable $\bm{g}(\bm{x})$ immediately yields $\bm{b}(\bm{x})$. 

To demonstrate the estimation, we explain the estimation procedure for $\bm{b}(\bm{x})$. Here, we introduce the following representation for the Koopman generator matrix $L$:
\begin{align}
L = 
\begin{bmatrix}
\ell_{11} & \cdots & \ell_{1 N_{\mathrm{dic}}}\\
\vdots & \ddots & \vdots\\
\ell_{N_{\mathrm{dic}}1} & \cdots & \ell_{N_{\mathrm{dic}} N_{\mathrm{dic}}}\\
\end{bmatrix}
=
\begin{bmatrix}
\bm{\ell}_1^{\top} \\
\vdots \\
\bm{\ell}_{N_{\mathrm{dic}}}^{\top} 
\end{bmatrix}.
\end{align}
Note that the matrix $L$ gives an approximation of the time derivative:
\begin{align}
(\mathcal{L}\bm{\psi})(\bm{x}) = \frac{d}{dt}\bm{\psi}\big(\bm{x}(t)\big) \simeq L \bm{\psi}\big(\bm{x}(t)\big). 
\label{eq_connect_two_ways}
\end{align}
In the following, we calculate this time derivative in two different ways.

First, we consider a simple application of $\mathcal{L}$. As exemplified in Eq.~\eqref{eq_dictionary_example}, the monomial dictionary has the constant $1$ as the first element. Then, the elements corresponding to the full-state observable are $\psi_{2}(\bm{x}) = x_1, \dots, \psi_{D+1}(\bm{x}) = x_D$ when we consider a $D$-dimensional state space. A naive application of $\mathcal{L}$ to $\psi_{i+1}(\bm{x})$ ($i = 1, \dots, D$) leads to
\begin{align}
(\mathcal{L} \psi_{i+1})(\bm{x})
=\sum_{d=1}^{D} b_d(\bm{x})\frac{\partial}{\partial x_d} x_{i}
=b_i(\bm{x}).
\label{eq_appcoximation_with_L_mid_1}
\end{align}
Then, we obtain
\begin{align}
(\mathcal{L}\bm{g})(\bm{x})
= \begin{bmatrix}
g_1(\bm{x})\\
\vdots\\
g_{D}(\bm{x})
\end{bmatrix}
= \begin{bmatrix}
\psi_{2}(\bm{x})\\
\vdots\\
\psi_{D+1}(\bm{x})
\end{bmatrix}
 = \begin{bmatrix}
b_1(\bm{x})\\
\vdots\\
b_D(\bm{x})
\end{bmatrix}.
\end{align}
Second, the right-hand side of Eq.~\eqref{eq_connect_two_ways} gives 
\begin{align}
\frac{d}{dt}\psi_{i+1}(\bm{x})
\simeq \bm{\ell}_{i+1}^{\top} \bm{\psi}(\bm{x}) .
\label{eq_approximation_with_L_row}
\end{align}
Hence, we finally obtain the following result:
\begin{align}
b_i(\bm{x}) = \bm{\ell}_{i+1}^{\top} \bm{\psi}(\bm{x}).
\label{eq_for_b_final}
\end{align}
Since we use the monomial dictionary functions, the interpretation of Eq.~\eqref{eq_for_b_final} is straightforward. For example, when we have
\begin{align}
\bm{\ell}^{\top}_1 = 
\begin{bmatrix}
0 & 4 & 0 & 0 & 0 & 0 & -4 & 0 & 0 & \cdots
\end{bmatrix},\\
\bm{\ell}^{\top}_2 = 
\begin{bmatrix}
0 & 0 & -2 & 0 & 0 & 0 & 0 & 0 & 0 & \cdots
\end{bmatrix},
\end{align}
the dictionary $\bm{\psi}(\bm{x})$ in Eq.~\eqref{eq_dictionary_example} yields
\begin{align}
b_1(\bm{x}) = 4 x_1 -4 x_1^3, \quad 
b_2(\bm{x}) = -2 x_2.
\end{align}

For the diffusion coefficient functions, a similar discussion is possible. Let the index $k$ denote the dictionary function $x_{i} x_{j}$, i.e., $\psi_k(\bm{x}) = x_{i} x_{j}$. Then, a naive application of $\mathcal{L}$ to $\psi_k(\bm{x})$ yields
\begin{align}
(\mathcal{L}\psi_k)(\bm{x})
&= \sum_{d=1}^{D} b_{d}(\bm{x})\frac{\partial}{\partial x_{d}}x_ix_j
+
\frac{1}{2}\sum_{l=1}^{D}\sum_{m=1}^{D}a_{lm}(\mathbf{x})
\frac{\partial^2}{\partial x_{l}\partial x_{m}} x_ix_k \nonumber \\
&= x_i b_j(\bm{x})+x_j b_i(\bm{x}) + a_{ij}(\bm{x}).
\label{eq_for_a}
\end{align}
The left-hand side of Eq.~\eqref{eq_for_a} is similarly approximated with Eqs.~\eqref{eq_appcoximation_with_L_mid_1} and \eqref{eq_for_b_final}. Then, we have
\begin{align}
a_{ij}(\bm{x})
&= (\mathcal{L}\psi_k)(\bm{x}) - x_j b_i(\bm{x}) - x_i b_j(\bm{x}) \nonumber\\
&\simeq \bm{\ell}_{k}^{\top}\bm{\psi}(\bm{x})
-x_j\left( \bm{\ell}_{i+1}^{\top}\bm{\psi}(\bm{x}) \right)
-x_i\left( \bm{\ell}_{j+1}^{\top}\bm{\psi}(\bm{x}) \right),
\label{eq_for_a_final}
\end{align}
where Eq.~\eqref{eq_for_b_final} is used to replace $b_i(\bm{x})$ and $b_j(\bm{x})$.

In our two-dimensional example, the concrete expression with the dictionary of Eq.~\eqref{eq_dictionary_example} is summarized as
\begin{align}
\bm{b}(\bm{x})=
\begin{bmatrix}
b_1(\bm{x})\\
b_2(\bm{x})
\end{bmatrix}
\simeq
\begin{bmatrix}
\bm{\ell}_2^{\top} \bm{\psi}(\bm{x})\\
\bm{\ell}_3^{\top} \bm{\psi}(\bm{x})
\end{bmatrix},
\end{align}
and 
\begin{align}
a_{11} &\simeq \bm{\ell}_4^{\top}\bm{\psi}(\bm{x}) - 2x_1 \big(\bm{\ell}_2^{\top}\bm{\psi}(\bm{x}) \big), \\
a_{12} = a_{21} &\simeq \bm{\ell}_5\bm{\psi}(\bm{x}) - x_2 \big(\bm{\ell}_2^{\top} \bm{\psi}(\bm{x}) \big) - x_1 \big(\bm{\ell}_3^{\top}\bm{\psi}(\bm{x})\big), \\
a_{22} &\simeq \bm{\ell}_6^{\top}\bm{\psi}(\bm{x}) - 2x_2 \big(\bm{\ell}_3^{\top}\bm{\psi}(\bm{x}) \big).
\label{eq:ax}
\end{align}

\subsection{Experimental settings in Ref.~\cite{Klus2020}}

In the numerical demonstration in Ref.~\cite{Klus2020}, it is assumed that $\{\bm{b}(\bm{x}_n)\}_{n=1}^N$ and $\{\Sigma(\bm{x}_n)\}_{n=1}^N$ (and hence, $\{A(\bm{x}_n)\}_{n=1}^{N}$) are known for a given set of training data points $\{\bm{x}_n\}_{n=1}^N$. Note that $\{\bm{b}(\bm{x}_n)\}_{n=1}^{N}$ and $\{A(\bm{x}_n)\}_{n=1}^{N}$ are scalar values that collect the values of time derivatives for each coordinate $\{\bm{x}_n\}_{n=1}^{N}$ of the data. Hence, the information does not directly yield the ``function'' $\bm{b}(\bm{x})$ for the time derivatives of all states $\bm{x}$. Then, the problem setting is still suitable for estimating the underlying system equations from data.

\begin{figure}
    \begin{center}
        \includegraphics[width=75mm]{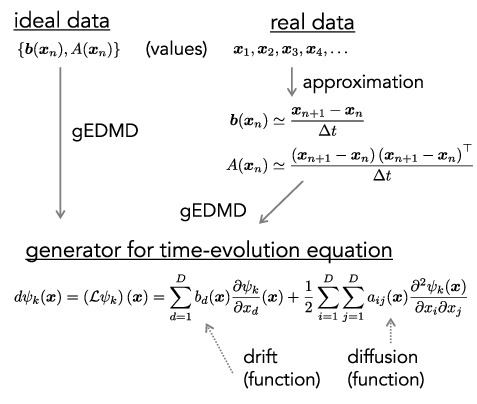}
        \caption{Framework of estimation. In an ideal situation, the data $\{\bm{b}(\bm{x}_n)\}_{n=1}^N$ and $\{A(\bm{x}_n)\}_{n=1}^N$ are assumed to be known. In a realistic case, we must estimate these quantities only from sample trajectories.}
        \label{fig_framework}
    \end{center}
\end{figure}

However, it would not be realistic to have the information of such data. Here, we consider an application of the gEDMD algorithm to sample trajectories, and hence $\{\bm{b}(\bm{x}_n)\}_{n=1}^{N}$ and $\{A(\bm{x}_n)\}_{n=1}^{N}$ are unknown. Actually, Ref.~\cite{Klus2020} also gave the discussion for such cases, where it was suggested to approximate $\{\bm{b}(\bm{x}_n)\}_{n=1}^{N}$ and $\{A(\bm{x}_n)\}_{n=1}^{N}$ from the time-series data $\{\bm{x}_{n}\}_{n=1}^{N}$ as
\begin{align}
\bm{b}(\bm{x}_n) &\simeq \frac{1}{\Delta t}(\bm{x}_{n+1}-\bm{x}_{n}),
\label{eq_b_from_naive_diff}\\
A(\bm{x}_n) &\simeq \frac{1}{\Delta t}(\bm{x}_{n+1}-\bm{x}_{n})(\bm{x}_{n+1}-\bm{x}_{n})^{\top},
\label{eq_a_from_naive_diff}
\end{align}
where $\Delta t$ is the time interval of the time-series data. That is, we employed the finite differences as the approximation. Figure~\ref{fig_framework} shows the framework of estimation. When the data $\{\bm{b}(\bm{x}_n)\}_{n=1}^N$ and $\{A(\bm{x}_n)\}_{n=1}^N$ are known, considerably accurate estimation is possible even in a small dataset; see Ref.~\cite{Klus2020}. However, in a realistic situation, these quantities must be estimated from sample trajectories. Since there was no demonstration of this approximation in Ref.~\cite{Klus2020}, we here check the ability of the approximation in preliminary numerical experiments.

\begin{figure}
    \begin{center}
        \includegraphics[width=75mm]{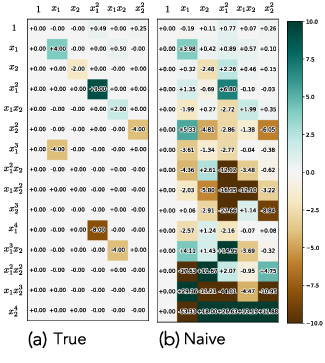}
        \caption{(Color online) The true matrix $L^{\top}$ and that obtained by the naive finite difference method in Eqs.~\eqref{eq_b_from_naive_diff} and \eqref{eq_a_from_naive_diff}. The details of the experimental settings are the same as in Sect.~5; we will denote them later.}
        \label{fig_result_naive}
    \end{center}
\end{figure}

Note that the data includes significant noise because of the stochastic nature, and hence the naive approach based on Eqs.~\eqref{eq_b_from_naive_diff} and \eqref{eq_a_from_naive_diff} is difficult to apply. Figure~\ref{fig_result_naive} shows the result of estimation for $L^{\top}$. The details of experimental settings are the same as those in Sect.~5; we will explain them later. Compared with the true answer, it is clear that many values are left where they should be zero. This is the effect of noise in the data. Of course, one could easily imagine using the sparse estimation, i.e., the least absolute shrinkage and selection operator (lasso) \cite{Tibshirani1996}. Then, the cost function in Eq.~\eqref{eq_cost} is replaced with
\begin{align}
L&=\argmin_{\widetilde{L}}
\left(\sum_{n=1}^{N}\left\|\bm{d\psi}(\bm{x}_n)-\widetilde{L}\bm{\psi}(\bm{x}_n)\right\|^2 + \lambda \left( \sum_{ij} \lvert \widetilde{\ell}_{ij} \rvert\right)
\right),
\label{eq_cost_lasso}
\end{align}
where $\lambda$ is a hyperparameter. We expect the regularization term with $\lambda$ enhances to increase the number of elements with a value of zero in the estimation. However, a naive application of the sparse estimation is time-consuming. In the numerical experiments, we use the totally $10$-th order monomials as the dictionary, i.e., $x_1^n x_2^m$ with $n+m \leq 10$. The dictionary size is $N_{\mathrm{dic}} = 66$, and in total, the Koopman generator matrix $L$ has $4,356$ elements. Hence, it is preferable to reduce the data size; here, we use $2\times 10^6$ data points. We will show the corresponding numerical results later.

A slightly different approach would be possible; we apply some preprocessing of the data to reduce the noise effects before applying the gEDMD algorithm. That is, instead of the naive finite difference method in Eqs.~\eqref{eq_b_from_naive_diff} and \eqref{eq_a_from_naive_diff}, the following conditional expectations are used:
\begin{align}
\bm{b}(\bm{x}) &\coloneqq\lim_{\Delta t \to 0}\mathbb{E}\left[\left.\frac{1}{\Delta t}(\bm{X}(t+\Delta t)-\bm{x})\right|\bm{X}(t)=\bm{x}\right], \\
A(\bm{x}) &\coloneqq\lim_{\Delta t \to 0}\mathbb{E}\left[\left.\frac{1}{\Delta t}(\bm{X}(t+\Delta t)-\bm{x})(\bm{X}(t+\Delta t)-\bm{x})^{\top}\right|\bm{X}(t)=\bm{x}\right].
\label{eq:kramers-moyal}
\end{align}
It is possible to estimate the conditional expectations from the dataset. Although Ref.~\cite{Klus2020} remarked on the conditional expectations, there are no comments on the practical implementation of the estimation and no numerical experiments with this formulation. As discussed later, the naive approach alone will not yield good results. In the following sections, we propose the method based on the conditional expectations.

\section{Proposal for Weighted Expectations via Data Classification by Clustering}

\subsection{Proposal 1: Clustering and weighted expectations}

We first consider a practical implementation to calculate the conditional expectations from data. It would be natural to employ the following locally weighted expectations:
\begin{align}
\widetilde{\bm{b}}(\bm{x})&=\frac{1}{\sum_{n=1}^{N-1} K_{H}(\bm{x},\bm{x}_n)}\sum_{n=1}^{N-1} K_{H}(\bm{x},\bm{x}_n)\left[\frac{1}{\Delta t}(\bm{x}_{n+1}-\bm{x}_n)\right],
\label{eq_weighted_bx}\\
\widetilde{A}(\bm{x})&=\frac{1}{\sum_{n=1}^{N-1} K_{H}(\bm{x},\bm{x}_n)} \nonumber \\
&\quad \times \sum_{n=1}^{N-1} K_{H}(\bm{x},\bm{x}_n)\left[\frac{1}{\Delta t}(\bm{x}_{n+1}-\bm{x}_n)(\bm{x}_{n+1}-\bm{x}_n)^{\top}\right],
\label{eq_weighted_ax}
\end{align}
where $K_{H}(\bm{x},\bm{x}_n)$ is a kernel function that yields the weighted averages. Here, we use the Gaussian kernel
\begin{align}
K_H(\bm{x},\bm{x}_n)=\exp \left(-\frac{1}{2}(\bm{x}_n-\bm{x})^\mathrm{T}H^{-1}(\bm{x}_n-\bm{x}) \right),
\label{eq:kh}
\end{align}
where $H$ represents a bandwidth matrix of the Gaussian kernel. As Proposal~1, we simply set $H$ as a diagonal matrix $\mathrm{diag}(h,\dots,h)$, where $h$ is a bandwidth parameter.

\begin{figure}
    \begin{center}
        \includegraphics[width=75mm]{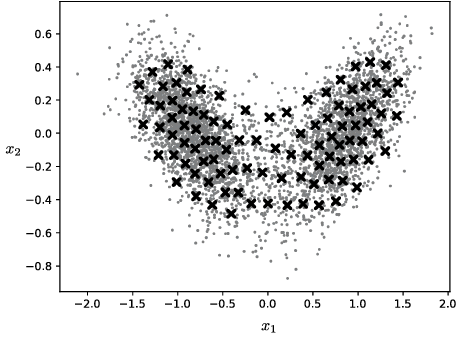}
        \caption{Selection of representative points by clustering. After excluding 5\% of points in the edge region, we select $100$ representative points; only the conditional expectations on these representative points are used in the gEDMD algorithm.} 
        \label{fig_representative_points}
    \end{center}
\end{figure}

When we apply the naive finite difference method in Eqs.~\eqref{eq_b_from_naive_diff} and \eqref{eq_a_from_naive_diff}, the noise in the data is reduced by the least squares method in the gEDMD. Then, we require a large amount of data. In contrast, in our proposed method, the noise in the data is reduced in the preprocessing using the weighted expectations. After the preprocessing, the gEDMD algorithm does not need so much data; there is no need to use all the data points $\{\bm{x}_{n}\}_{n=1}^{N}$ for the gEDMD algorithm. Here, we use only a small number of representative points, and the weighted expectations $\widetilde{\bm{b}}(\bm{x}), \widetilde{A}(\bm{x})$ are calculated only for these representative points.

There are various methods to select the representative points, such as grid points and randomly sampled points from the data. Here, we select the representative points by clustering with the $k$-means algorithm \cite{Hastie_book}. With the selection of representative points using clustering, the major parts of the space can be covered. Although defining grid points as representative points also covers the major parts of the space, it may fail to estimate weighted expectations accurately if there are insufficient points around the defined grid points. Moreover, while randomly sampling points from the data to use as representative points may lead to a concentration of representative points, the clustering-based selection prevents such concentration.

Figure~\ref{fig_representative_points} shows an example of the representative points. The dots represent original data points, whereas the crosses represent the representative points. In the $k$-means algorithm, it is necessary to specify the number of clusters beforehand, and we set the number as that of the representative points. In Fig.~\ref{fig_representative_points}, we set the number as $100$. In addition, original data points in the edge region may not get weighted expectations well; there may not be enough nearby data points. Hence, the outlier detection method, \verb|IsolationForest| in Ref.~\cite{Liu2008}, is employed to exclude the 5\% of points in the edge region.

\subsection{Proposal 2: Additional clustering and covariance estimation}

There is a further proposal that plays a crucial role in the numerical experiments. One could expect that Proposal~1 will reduce the effect of noise. However, we sometimes require slightly careful calculations for the weighted expectation, especially in the case of state-dependent diffusion terms. As we saw in Fig.~\ref{fig_sample}, there are some state-dependent features in the sample trajectory; some parts of the structure are tilted at an angle to the right, while others are tilted to the left. It would be better to take the weighted expectations using only data points undergoing similar motions. Hence, apart from clustering to select representative points, we perform an additional clustering process.

\begin{figure}[tb]
    \begin{center}
        \includegraphics[width=75mm]{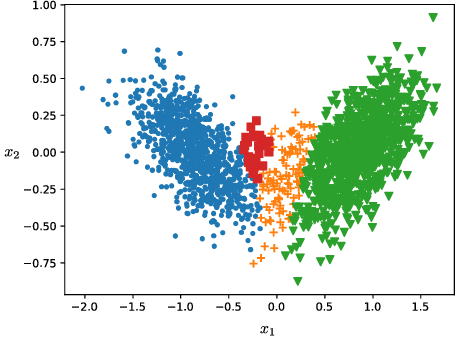}
        \caption{(Color online) Clustering for the calculation of covariance estimation. The estimated covariance is used as the bandwidth in the locally weighted expectations.}
        \label{fig_dpmm}
    \end{center}
\end{figure}

To address this issue, we employ the Dirichlet process mixture model (DPMM) \cite{Ferguson1973} to obtain rough clustering results of the data with the same tendency. The DPMM clustering can appropriately classify the data even when the number of clusters is unknown. Figure~\ref{fig_dpmm} shows an example of the clustering results, in which we finally have four clusters. The DPMM yields the appropriate clustering number automatically; this characteristic is preferable when the data dimension is high and we cannot visualize the data points. The covariance matrix for each cluster is evaluated and used as $H$ in the Gaussian kernels.

After the clustering with the DPMM, we perform a procedure similar to that in Proposal~1; the representative points are selected by the $k$-means method. For each representative point, we determine which cluster it belongs to in the DPMM and calculate the weighted expectation using the data points that belong to the same cluster. We set the bandwidth matrix $H$ with the covariance matrix of the cluster to which the representative point belongs.

\section{Numerical Experiments}

\subsection{Settings}

Time-series data for the double-well potential model is generated through simulations of stochastic differential equations. The Euler--Maruyama method is used in the simulations \cite{Kloeden_book}, in which the time discretization is $\Delta t = 1.0\times 10^{-3}$. All data points are used as the sample trajectory, and we finish the time evolution with $T=2000$. Then, the number of data points is $N=2\times10^{6}$.

The total maximum degree of monomial functions in the dictionary $\bm{\psi}(\bm{x})$ is $10$; i.e., $x_1^n x_2^m$ with $n+m \leq 10$. Then, there are a total of $N_{\mathrm{dic}}=66$ dictionary functions.

The details of the settings for each method are as follows.
\begin{itemize}
\item The method of ``Naive lasso'' uses the naive finite difference in Eqs.~\eqref{eq_b_from_naive_diff} and \eqref{eq_a_from_naive_diff}. Since all the $2\times10^{6}$ data points are too large for the lasso, we select every 100 pieces of the original data. Then, we analyze the selected $N=2\times10^{4}$ data points; the gEDMD algorithm finishes in about $2$ min in a MacBookAir with an M2 processor. The hyperparameter for the lasso is set to $\lambda = 0.01$; we tried some parameters, and we show the results of the parameters that give moderate results.
\item The method of ``Proposal~1'' is described in Sect.~4.1. The procedure of representative point selection is relatively time-consuming. Hence, we select every $100$ pieces of the original data and apply the $k$-means clustering to the $2\times10^{4}$ data points. Note that we take the weighted expectations with all data points with $N=2\times10^{6}$. After evaluating the weighted expectations on the representative points, we apply the gEDMD algorithm to only the $100$ representative points; the total calculation time is less than $10$ s. The isotropic bandwidth is $h=0.2$, and the hyperparameter for the lasso is $\lambda = 1.0\times 10^{-6}$.
\item The method of ``Proposal~2'' is described in Sect.~4.2. As for the DPMM algorithm, \verb|BayesianGaussianMixture| in \verb|scikit-learn| \cite{scikit_learn} is used. The number of mixture components is set as $10$, and finally, four components are selected in our dataset, as shown in Fig~\ref{fig_dpmm}. Other settings are the same as those in ``Proposal~1.'' The total computation time is about $40$ s, in which the DPMM takes about $30$ s. The hyperparameter for the lasso is $\lambda = 1.0\times 10^{-3}$.
\end{itemize}

\subsection{Numerical results}

\begin{figure*}
    \begin{center}
        \includegraphics[width=140mm]{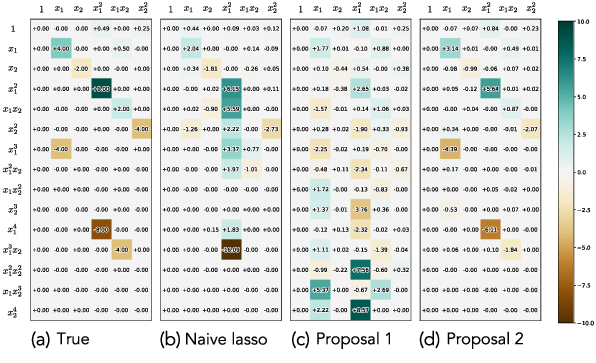}
        \caption{(Color online) Koopman generator matrix $L^{\top}$. Only elements with small total degrees are shown. (a) shows the true values of the model. (b), (c), and (d) correspond to results obtained from ``Naive lasso,'' ``Proposal~1,'' and ``Proposal~2,'' respectively.}
        \label{fig_results_matrix}
    \end{center}
\end{figure*}

\begin{figure*}
    \begin{center}
        \includegraphics[width=140mm]{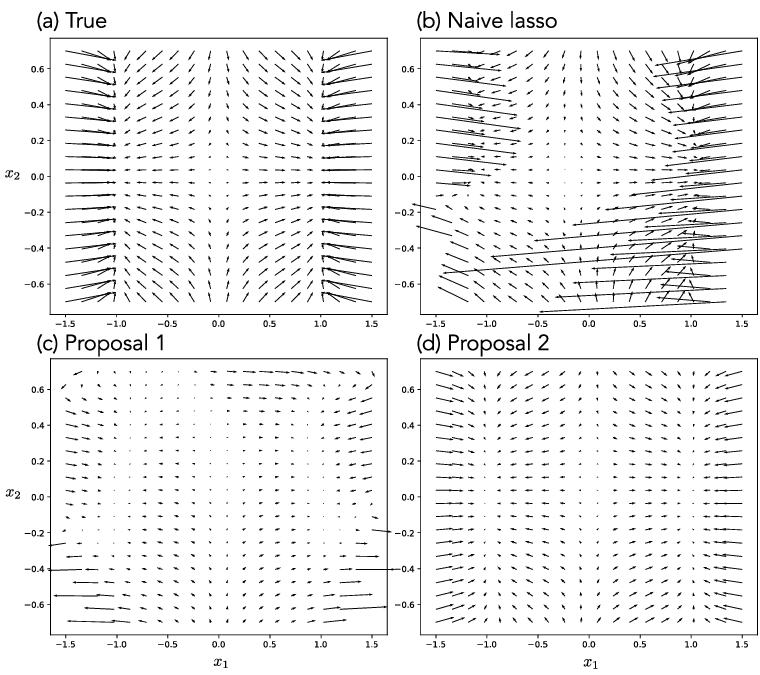}
        \caption{Evaluated drift coefficient functions $\bm{b}(\bm{x})$ on certain grid points. (a) is evaluated using Eq.~\eqref{eq_drift}. (b), (c), and (d) correspond to results obtained from ``Naive lasso,'' ``Proposal~1,'' and ``Proposal~2,'' respectively.}
        \label{fig_results_b}
    \end{center}
\end{figure*}

Figure~\ref{fig_results_matrix}(a) shows a portion of the true Koopman generator matrix. Figures~\ref{fig_results_matrix}(b)--7(d) are the Koopman generator matrices obtained from ``Naive lasso,'' ``Proposal~1,'' and ``Proposal~2,'' respectively. We see that there are still values left for elements that should have a value of zero. Note that the true value in row $x_1$ and column $x_1$ is $+4$, although (b)--(d) yield smaller values than $+4$. These results indicate that sparsity is sufficient. If we increase the hyperparameter of sparsity, the number of elements with zero value will increase. However, the element in row $x_1$ and column $x_1$ becomes far from the true value. From Fig.~\ref{fig_results_matrix}, it is clear that the additional clustering process in ``Proposal~2'' has significantly improved the results.

Figure~\ref{fig_results_b} shows more clearly the effect of the additional clustering process in ``Proposal~2,'' which depicts the evaluated drift coefficient functions $\bm{b}(\bm{x})$ on certain grid points. In the regions near $x_1 \simeq -1.5$ and $x_2 \simeq 1.5$, ``Naive lasso'' and ``Proposal~1'' yield some arrows with strange directions and lengths. In contrast to these results, ``Proposal~2'' shows a behavior similar to the true one. The arrows are short in the regions near $x_1 \simeq -1.5$ and $x_2 \simeq 1.5$; we believe the reason is that there are insufficient data points around these regions; see Fig.~\ref{fig_sample}.

\begin{figure}
    \begin{center}
        \includegraphics[width=70mm]{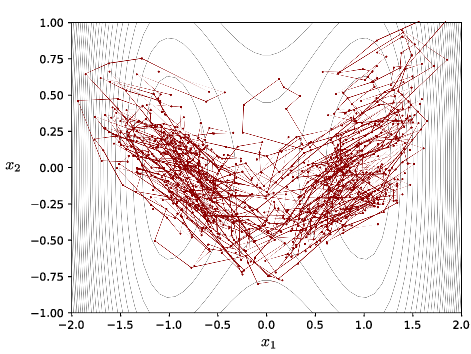}
        \caption{(Color online) A sample trajectory obtained from the estimated model in ``Proposal~2.''}
        \label{fig_sample_trajectory_model}
    \end{center}
\end{figure}

Figure~\ref{fig_sample_trajectory_model} shows a sample trajectory obtained from the estimated model in ``Proposal~2.'' Although the trajectory is a slightly larger than the original one, it shows the going back-and-forth behavior within the two valleys of the potential. It also appropriately reproduces the deviation from the potential shape due to the state dependence of the diffusion term.

The example here yields the sparse Koopman matrix shown in Fig.~\ref{fig_results_matrix}(a). We show numerical experimental results for a slightly different example in the Appendix.

\section{Conclusion}

One may consider that the naive weighted expectations with the surrounding data points would improve the performance. However, the numerical experiments in this study show that one more idea yields a performance improvement; the local information extracted via the additional clustering is crucial in the Koopman generator estimation using only the sample trajectories. There are many works on the EDMD and gEDMD algorithms, including analysis using eigenvalue modes and more accurate prediction with the aid of deep neural networks. However, the research on model characterization at the equation level, which is of interest from a physics perspective, has only just begun. We believe that the technique proposed here will help expand the range of future applications of the gEDMD algorithm.

There are some remaining issues for future studies of model estimation. One of them is the choice of the bandwidth matrix. In this study, we employed the covariance matrix as the bandwidth matrix. One may consider that the bandwidth matrix selection method developed in the kernel density estimation \cite{Silverman_book} may be applicable instead of the proposed simple method. However, in our preliminary numerical experiments, the bandwidth matrix selection method alone did not work well; we require additional parameters to obtain good estimation results. Hence, the development of a bandwidth matrix selection method for the weighted expectations with theoretical support may be promising. Another problem is the curse of dimensionality; the size of the dictionary becomes enormous as the dimension of the space increases. There are several studies that address this issue using the tensor train format \cite{Klus2018,Nuske2021,Lucke2022}, and the use should be considered in the future.

\acknowledgements
This work was supported by JSPS KAKENHI Grant Number JP23H04800.

\appendix

\section{A slightly dense case}

The example in Sect.~5 yields the sparse Koopman matrix shown in Fig.~\ref{fig_results_matrix}. Here, we consider the following potential $V(\bm{x})$ instead of Eq.~\eqref{eq_potential}:
\begin{align}
V(\bm{x}) =& - 0.4 x_1 + 0.4 x_2 - x_1^2 - 0.3 x_1 x_2 + 2.0 x_2^2 \nonumber \\
&+ 0.2 x_1^3 + 0.4 x_1^2 x_2 - 0.4 x_1 x_2^2 - 0.2 x_2^3  \nonumber \\
&+ x_1^4 - 0.2 x_1^3 x_2 + 0.2 x_1^2 x_2^2 + 0.2 x_2^4,
\label{eq_appendix_example}
\end{align}
which yields a slightly dense Koopman matrix. In the estimation, we change the hyperparameter $\lambda$ as $\lambda = 1.0\times 10^{-8}$ for ``Proposal~1'' and $\lambda = 1.0\times 10^{-4}$ for ``Proposal~2.'' Other settings are the same as in Sect.~5.

\begin{figure*}
    \begin{center}
        \includegraphics[width=140mm]{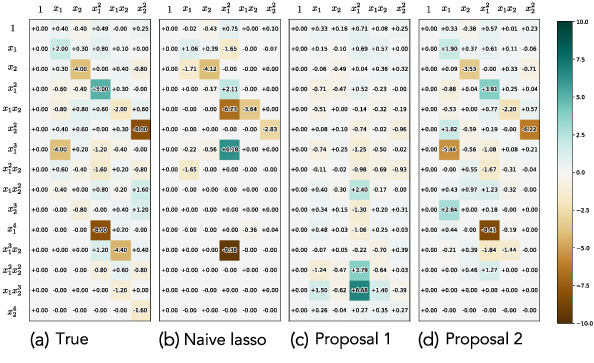}
        \caption{(Color online) Koopman generator matrix $L^{\top}$ for the case with Eq.~\eqref{eq_appendix_example}. Only elements with small total degrees are shown. (a) shows the true values of the model. (b), (c), and (d) correspond to results obtained from ``Naive lasso,'' ``Proposal~1,'' and ``Proposal~2,'' respectively.}
        \label{fig_results_matrix_appendix}
    \end{center}
\end{figure*}

The numerical results are shown in Fig.~\ref{fig_results_matrix_appendix}. We see that ``Proposal~2'' yields reasonable results compared with other methods, even in the slightly dense example.

\begin{acknowledgments}
This work was supported by JSPS KAKENHI Grant Number JP23H04800.
\end{acknowledgments}

\end{document}